# Conformal Killing $L^2$-forms on complete Riemannian manifolds with nonpositive curvature operator

Sergey Stepanov and Irina Tsyganok

We give a classification for connected complete locally irreducible Riemannian manifolds with nonpositive curvature operator, which admit a nonzero closed or co-closed conformal Killing $L^2$-form. Moreover, we prove vanishing theorems for closed and co-closed conformal Killing $L^2$-forms on some complete Riemannian manifolds.

## 1. Introduction and results

*Conformal Killing forms* (also called as conformal Killing-Yano tensors) have been defined on Riemannian manifolds more than forty-five years ago by S. Tachibana and T. Kashiwada (see [24] and [9]) as a natural generalization of conformal Killing vector fields. We also know from the literature about *closed conformal Killing forms* or, otherwise, closed conformal Killing-Yano tensors and *co-closed conformal Killing forms* or, otherwise, Killing-Yano tensors (see, for example, [7, pp. 426-427]; [18; pp. 559-564]). We remark here that the Hodge dual of a co-closed Killing form is a closed conformal Killing form. Moreover, the converse is also true (see [14]). Surveys of the publications on conformal Killing, co-closed and closed conformal Killing forms and their numerous applications can be found in the introductions to our papers [19] and [20]. In addition, it should be taken into account the list of recent papers on these forms: [4]; [6]; [13]; [21]; [22] and [27]. In the present paper we consider conformal Killing, co-closed and closed conformal Killing $L^2$-forms of degree $p$ for $1 \leq p \leq n-1$ on a simply connected and complete Riemannian manifold $(M, g)$ of dimension $n$ for $n \geq 2$.

Here we regard the Riemannian curvature tensor $Rm$ of $(M,g)$ as a symmetric algebraic operator $\bar{R} : \Lambda^2(T_xM) \to \Lambda^2(T_xM)$ on the vector space $\Lambda^2(T_xM)$ of 2-forms over tangent space $T_xM$ at an arbitrary $x \in M$ (see [11, pp. 36-37]). We say that the manifold $(M, g)$ has a nonpositive (respectively nonnegative) curvature operator $\bar{R}$ if $g(\bar{R}(\theta),\theta) \leq 0$ (respectively $g(\bar{R}(\theta),\theta) \geq 0$) for all two-forms $\theta \neq 0$. There have been many papers on the relationship between the curvature operator $\bar{R}$ of a Riemannian manifold $(M,g)$ and some global characterization of it, such as its homotopy type, topological types and etc.

In connection with above, the first our result on conformal Killing forms will be proved by the most important analytic method of differential geometry "in the large" which derived by S. Bochner for proving so-called vanishing theorems under appropriate curvature conditions on compact Riemannian manifolds (see [28]). S.-T. Yau generalized this method of proving vanishing theorems for the case of complete noncompact Riemannian manifolds (see, for example, [29]). We use a generalization of the "Bochner technique" to prove the following statement.

**Lemma 1.** *Suppose the curvature operator $\bar{R}$ is negative semi-define at every point of a complete non-compact Riemannian manifold $(M,g)$. Then every closed or co-closed conformal Killing $L^2$-form on $(M,g)$ is a parallel form. If either the volume of $(M,g)$ is infinite, then every closed or co-closed conformal Killing $L^2$-form on $(M,g)$ is identically zero.*

If $\omega$ is a parallel form on $(M,g)$ then $\omega_x$ is invariant under the holonomy representation $\mathrm{Hol}(x)$ at each point $x \in M$. Therefore, finding on a Riemannian manifold $(M,g)$ parallel forms is equivalent to the algebraic problem of finding the invariant forms of the holonomy representation $\mathrm{Hol}(g)$. We give an example of the implementation of this principle, which is the main result of the present paper.

**Theorem 1**. *Let a simply connected complete non-compact Riemannian manifold $(M,g)$ with nonpositive curvature operator admits a non-zero not decomposable closed or co-closed conformal Killing $L^2$-form. If in addition, we assume that $(M,g)$ is a locally irreducible Riemannian manifold, then it is a manifold of the following list: globally symmetric manifolds of non-compact type, Kähler and quaternionic-Kähler manifolds.*

**Remark 1**. In addition, we recall that an arbitrary co-closed conformal Killing $p$-form $\omega$ on an $n$-dimension compact Kählerian manifold is parallel for $2 \leq p \leq n-1$ (see [25]). From this statement follows that an arbitrary closed conformal Killing $(n-p)$-form $\omega$ on an $n$-dimension compact Kählerian manifold is parallel as well.

In particular, for $\dim M = 2p$ we prove the following lemma using a generalization of the "Bochner technique".

**Lemma 2.** *Suppose the curvature operator $\bar{R}$ is negative semi-define at every point of a $2p$-dimensional complete non-compact Riemannian manifold $(M,g)$. Then every closed or co-closed conformal Killing $L^2$-form of degree p is parallel form on $(M,g)$. If either the volume of $(M,g)$ is infinite, then every closed or co-closed conformal Killing $L^2$-form of degree p is identically zero on $(M,g)$.*

**Remark 2**. In turn, we have proved that an arbitrary conformal Killing $p$-form $(1 \leq p \leq n-1)$ is parallel on an $n$-dimension compact Riemannian manifold with nonpositive curvature operator (see [17]).

The following theorem is an analogue of Theorem 1.

**Theorem 2**. *Let a $2p$-dimensional simply connected complete non-compact Riemannian manifold $(M,g)$ with nonpositive curvature operator admits a non-zero not decomposable conformal Killing $L^2$-form of degree p. If in addition, we assume that $(M,g)$ is a locally irreducible, then it is a manifold of the following list: globally symmetric manifolds of non-compact type, Kähler and quaternionic-Kähler manifolds.*

We can formulate the converse statement in one particular case.

**Corollary 1.** *A globally symmetric space of non-compact type $(M,g)$ with infinite volume does not admit a nonzero conformal Killing $L^2$-form.*

**Remark 3**. In addition, we know that if a symmetric space $(M,g)$ of *compact type* carries a non-parallel co-closed conformal Killing $p$-form for $p \geq 2$, then its universal cover $\tilde{M}$ is either a round sphere, or has a factor isometric to a round sphere in its de Rham decomposition (see [12]).

An $n$-dimensional $(n \geq 2)$ Riemannian manifold $(M, g)$ is *generic* (see [1, p. 291]) if each point $x \in M$ admits the orthonormal base $\{e_i\}$ of the tangent space $T_x M$ such that the $n(n-1)$ endomorphisms $R(e_i, e_j)$ for the curvature tensor $R$ of $(M, g)$ will be linearly independent in the Lie algebra of the orthogonal group $O(n)$. In this case, the corollary is true.

**Corollary 2**. *A generic complete Riemannian manifold with nonpositive curvature operator does not admit nonzero closed and co-closed conformal Killing $L^2$-forms.*

In [1, p. 419] was proved the statement on Kähler manifold with parallel form that is not proportional to any degree of the Kähler 2-form. Here we interpret this statement as the following

**Corollary 3**. *Let a complete non-compact locally irreducible Kähler manifold admits a nonzero irreducible closed or co-closed conformal Killing $L^2$-form $\omega$. If the curvature operator $\bar{R} \leq 0$ everywhere on the Kähler manifold and $\omega$ is not some power of the Kähler form then the Kähler manifold is an Einstein manifold.*

## 2. Preliminary information

Let $\Lambda^p M$ be the bundle of differential $p$-forms over a connected complete Riemannian manifold $(M, g)$ of dimensions $n \geq 2$ and $C^\infty(\Lambda^p M)$ be the space of $C^\infty$-sections of the bundle $\Lambda^p M$ for $1 \leq p \leq n-1$. By $\nabla$ we will denote the Levi-Civita connection on $(M, g)$ and by $\nabla^*$ the canonical formal adjoint of $\nabla$. In addition, let $d: C^\infty(\Lambda^p M) \to C^\infty(\Lambda^{p+1} M)$ be the well known operator of exterior derivative and $d^*: C^\infty(\Lambda^p M) \to C^\infty(\Lambda^{p-1} M)$ be the codifferentiation operator which defined as the canonical formal adjoint of $d$. We recall that if $d\omega = 0$, then the $p$-form $\omega \in C^\infty(\Lambda^p M)$ is said to be *closed*. In tern, if $d^*\omega = 0$, then the $p$-form $\omega \in C^\infty(\Lambda^p M)$ is said to be *co-closed*. Using these operators, one constructs the well-known Hodge-de Rham Laplacian $\Delta = d^*d + dd^*$ which admits a Weitzenböck decomposition (see [1, p. 53]; [8, p. 144]; [11, p. 211])

(2.1) $$\Delta \omega = \nabla^*\nabla \omega + E(\omega)$$

for any $\omega \in C^\infty(\Lambda^p M)$ and an algebraic symmetric operator $E : \Lambda^p M \to \Lambda^p M$ that depends linearly in a known way (see, for example, [23]; [26]) on the Riemannian curvature tensor $Rm$ and the Ricci tensor $Ric$ of $(M,g)$. Further, a direct calculation based on (2.1) yields the well-known Bochner-Weitzenböck formula (see, for example, [26])

(2.2) $$\frac{1}{2}\Delta \|\omega\|^2 = g(\Delta \omega, \omega) - \|\nabla \omega\|^2 - g(E(\omega), \omega).$$

Let $\{e_1,...,e_n\}$ be locally defined orthonormal frame fields of the tangent bundle $TM$ and $\{\omega_1,...,\omega_n\}$ be its dual coframe fields. Then in terms of the basis, we can write $\omega = a_{i_1...i_p}\omega_{i_1} \wedge ... \wedge \omega_{i_p}$ where the summation is being performed over the multi-index $I = (i_1,...,i_p)$ for $1 \leq i_1,...,i_p \leq n$. With this understanding, we can write (see [10, p. 29-30])

(2.3) $$g(E(\omega), \omega) = g(\bar{R}(\bar\omega), \bar\omega)$$

for the 2-form $\bar\omega = a_{i_1...j_\alpha...i_p} \omega_{j_\alpha} \wedge \omega_{i_\alpha}$.

J.-P. Bourguignon in [2] proved the existence of the basis $\{d, d^*, D\}$ in the space of natural (with respect to isometric diffeomorphisms) first-order differential operators on the space $C^\infty(\Lambda^p M)$ with value in the space of homogeneous tensor on $(M,g)$. The third basis operator $D$ was not given explicitly. In our papers [14] and [15] we used some classical theorems of H. Weyl about the representation theory of the orthogonal group $O(n)$ and shown that the third operator $D$ has the form

$$D = \nabla - \frac{1}{p+1}d - \frac{1}{n-p+1}g \wedge d^*$$

where $(g \wedge d^*)(X_0, X_1, ..., X_p) = \sum_{a=1}^{p} (-1)^a g(X_0, X_a)(d^*\omega)(X_1, ..., X_{a-1}, X_{a+1}, ..., X_p)$

for any $\omega \in C^{\infty}(\Lambda^p M)$ and $X_1, ..., X_p \in C^{\infty}(TM)$. We also proved that the kernel of $D$ consists of conformal Killing $p$-forms. Therefore, a closed conformal Killing form $\omega$ is determined by the conditions $D\omega = 0$ and $d\omega = 0$. In turn, the conditions $D\omega = 0$ and $d^*\omega = 0$ define a co-closed conformal Killing form $\omega$. In addition, we proved in [14] that if $\omega$ is a closed conformal Killing $p$-form then $*\omega$ must be a co-closed conformal Killing $(n-p)$-form for the Hodge star operator $* : \Lambda^p M \to \Lambda^{n-p} M$, which is a familiar isomorphism of vector bundles $\Lambda^p M$ and $\Lambda^{n-p} M$ (see, for example, [1, p. 33]). Moreover, the converse is also true (see also [14]).

In [16] we found the operator $D^*$ formally adjoint to $D$ and then constructed the second-order differential operator

(2.4) $$D^*D = \frac{p}{p+1}\left(\nabla^*\nabla - \frac{1}{p+1}d^*d - \frac{1}{n-p+1}dd^*\right)$$

with $1 \leq p \leq n-1$. Properties of the operator $D^*D$ were studied in the following papers [16]; [19]; [20] and [21]. In particular, we proved in [16] that $D^*D$ is a second-order self-adjoint elliptic differential operator acting on $C^{\infty}(\Lambda^p M)$.

### 3. Proofs of the statements

A direct calculation yields the second inequality of Kato (see [2, p. 380])

(3.1) $$\|\omega\| \Delta \|\omega\| \leq g(\nabla^*\nabla\omega, \omega)$$

for an arbitrary form $\omega \in C^{\infty}(\Lambda^p M)$. By virtue of the Weitzenböck decomposition (2.1), the Kato inequality (3.1) can be written in the form

(3.2) $$\|\omega\| \Delta \|\omega\| \leq g(\Delta\omega, \omega) - g(E(\omega), \omega).$$

At the same time, by the Weitzenböck decomposition (2.1), the operator $D^*D$ can be written in the following form

(3.3) $$D^*D\omega = \frac{p}{p+1}\left(\frac{n-p}{n-p+1}\Delta\omega - E(\omega) - \frac{n-2p}{(p+1)(n-p+1)}d^*d\omega\right)$$

for an arbitrary $p$-form $\omega \in C^\infty(\Lambda^p M)$. Then from (3.3) we obtain equation

(3.4) $$\Delta\omega = \frac{n-p+1}{n-p}E(\omega)$$

for a closed conformal Killing $p$-form $\omega$. Using (3.4) we write the inequality (3.2) in the form

(3.5) $$\|\omega\|\Delta\|\omega\| \leq \frac{1}{n-p}g(E(\omega),\omega)$$

for a closed conformal Killing $p$-form $\omega$. If the curvature operator is nonpositive, then based on (2.3) we obtain from (3.5) the inequality

$$\|\omega\|(-\Delta\|\omega\|) \geq 0.$$

For a Riemannian manifold $(M,g)$ we have the "natural" Hilbert space $L^2(M) = L^2(M, d\operatorname{Vol}_g)$ where $\operatorname{Vol}_g$ is the Riemannian volume form associated to the metric $g$. In our case, assume that $\omega$ is a $L^2$-form on a complete $(M,g)$, then by Yau theorem (see [29, p. 664]) the function $\|\omega\|$ is constant for an arbitrary closed conformal Killing form $\omega$. Then it follows from (3.5) that $g(E(\omega),\omega) = 0$. In this case, we obtain from (3.4) that $g(\Delta\omega,\omega) = 0$. Finally, we obtain from the Weitzenböck formula (2.2) that $\nabla\omega = 0$ for an arbitrary closed conformal Killing $L^2$-form $\omega$ on a complete manifold $(M,g)$ with nonpositive curvature operator. Then $*\omega$ is a co-closed conformal Killing form whose covariant derivative vanishes. In particular, it means that $\|*\omega\|^2 = const$. In this case, we have $\operatorname{Vol}_g(M) < \infty$ for the volume $\operatorname{Vol}_g(M)$ of $(M,g)$ because $\int_M \|\omega\|^2 d\operatorname{Vol}_g < \infty$ for an arbitrary $L^2$-form $\omega$. Therefore, if the volume $\operatorname{Vol}_g(M)$ of a complete manifold $(M,g)$ with nonpositive curvature operator is infinite then $(M,g)$ does not admit a nonzero closed (resp. co-closed) conformal Killing $L^2$-form $\omega$. This completes the proof of Lemma 1.

We recall that the description of all the connected complete locally irreducible Riemannian manifolds with a form whose covariant derivative vanishes is given in Corollary 10.110 from the monograph [1, p. 306-307]. Namely, the following statement holds. Let $(M,g)$ be a connected complete locally irreducible Riemannian manifold admitting a nonzero and not decomposable form $\omega$ whose covariant derivative is vanishing. Then $(M,g)$ should fall in one of the following categories: locally symmetric manifolds, Kähler manifolds, quaternionic-Kähler manifolds, manifolds with holonomy group $G_2$ or $Spin(7)$. We recall that any complete, simply-connected locally symmetric space is globally symmetric. In addition, a Riemannian globally symmetric space of *non-compact type* has a nonpositive sectional curvature. Note that symmetric spaces of non-compact type are non-compact. Now the assertion about symmetric spaces in our theorem becomes obvious. Next, we note here that manifolds with holonomy group $G_2$ or $Spin(7)$ are Ricci-flat manifolds. On the other hand, every non-compact complete Riemannian manifold with non-negative Ricci curvature has infinite volume (see [29, p. 667]). Therefore, any non-compact complete Riemannian manifolds with holonomy group $G_2$ or $Spin(7)$ do not admit a not decomposable parallel $L^2$-form. These remarks complete the proof of our Theorem 1.

**Remark 4**. S.-T. Yau used the Laplacian $\Delta$ in the form $\Delta = trace_g \nabla^2$ without the minus sing (see [29]) that we took into account in our proof.

In particular, if $n = 2p$ we obtain from (3.4) that

$$(3.7) \qquad D^*D\omega = \frac{p}{p+1}\left(\frac{p}{p+1}\Delta\omega - E(\omega)\right).$$

Suppose that $\omega$ is a conformal Killing $p$-form, then it follows from (3.7) that $\Delta\omega = (p+1)p^{-1}E(\omega)$. Then based on (2.3) we obtain the inequality $\|\omega\|(-\Delta\|\omega\|) \geq 0$ if $\overline{R} \leq 0$ everywhere on $(M,g)$. Further, it can be shown that

$\nabla\omega=0$ if $\omega$ is a conformal Killing $L^2$-form on a complete manifold $(M,g)$ with nonpositive curvature operator. In this case, the assertions of Lemma 2 and Theorem 2 become obvious.

A Riemannian globally symmetric space $(M,g)$ is complete. We also know that a Riemannian symmetric space has nonpositive curvature operator if and only if it has nonpositive sectional curvature (see [5]). After the above remarks, the assertion of Corollary 1 becomes obvious.

Let $\text{Hol}^0(x)$ be the *restricted holonomy group* of $(M,g)$ at an arbitrary point $x \in M$ (see [1, p. 280]). Note that $\text{Hol}^0(x)$ is always contained in special orthogonal group $SO(n)$. But for a generic Riemannian manifold $\text{Hol}^0 = SO(n)$ at all its points (see [1, p. 291]). In this case, there is no a nontrivial parallel form $\omega$, i.e. there is no a nontrivial form $\omega$ with zero covariant derivative (see [1, p. 306]). This implies the validity of Corollary 2.

In [1, p. 307]) was proved the following statement: "Let $(M,g)$ be a Riemannian manifold admitting an exterior form $\omega$ whose covariant derivative is vanishing. Assume now that $\omega$ is not zero, not decomposable and not some power of the Kähler form if $(M,g)$ is Kähler. Then $(M,g)$ is automatically an Einstein manifold". Our Corollary 3 is a corollary of Theorem 1 and this statement.

## 4. Appendix

A form $\omega \in C^\infty(\Lambda^p M)$ is said to be *harmonic* if $\Delta\omega = 0$. According to Yau proposition (see [29, p. 663]) an $L^2$-form $\omega$ on a complete Riemannian manifold $(M,g)$ is harmonic if and only if it is closed and co-closed. The Bochner-Weitzenböck formula (2.2) implies that any harmonic $p$-form $\omega$ satisfies the inequality

$$(4.1) \qquad -\frac{1}{2}\Delta\|\omega\|^2 = \|\nabla\omega\|^2 + g(E(\omega),\omega) \geq 0.$$

on $(M,g)$ with the nonnegative curvature operator $\overline{R}$. Then by Yau theorem (see [29, p. 663]) the function $\|\omega\|$ is constant for an arbitrary harmonic form $\omega$ on a complete Riemannian manifold $(M,g)$ with $\overline{R} \geq 0$. An elementary linear algebra argument show that if the curvature operator in non-negative, then all the sectional curvatures of $(M,g)$ are also non-negative. Then the Ricci curvature of $(M,g)$ is non-negative as well. In this case, the volume $\text{Vol}_g(M)$ is infinite if $(M,g)$ is a non-compact complete manifold (see [29, p. 667]). As a result, a complete non-compact Riemannian manifold $(M,g)$ with non-negative curvature operator does not admit a nonzero harmonic $L^2$–form $\omega$. Therefore, we can not formulate a theorem on harmonic forms which is an analog of our Theorem 1. But the following statements hold.

**Corollary 4.** *A globally symmetric space of compact type $(M,g)$ does not admit a non-parallel harmonic form.*

**Proof.** It is well known that a Riemannian globally symmetric space of compact type $(M,g)$ has a nonpositive sectional curvature and also $(M,g)$ is compact. Then from (4.1) we obtain $\|\omega\| = const$ for an arbitrary harmonic form $\omega$ on a Riemannian globally symmetric space of compact type $(M,g)$. In this case, the covariant derivative of $\omega$ is vanishing.

**Corollary 5**. *A generic complete Riemannian manifold with non-negative curvature operator does not admit nonzero harmonic $L^2$-forms.*

**Proof.** It is well known that for a generic Riemannian manifold $\text{Hol}^0 = SO(n)$ at all its points (see [1, p. 291]). In this case, there is no a nontrivial form $\omega$ with zero covariant derivative (see [1, p. 306]). On the other hand, from (4.1) we obtain $\nabla \omega = 0$ for an arbitrary harmonic $L^2$-forms $\omega$ on a complete Riemannian manifold with non-negative curvature operator. This contradiction proves Corollary 5.

**Acknowledgments**

Our work was supported by RBRF grant 16-01-00053-a (Russia).

DEPARTAMENT OF MATHEMATICS FINANCE UNIVERSITY

UNDER THE GAVERMENT OF RUSSIAN FEDERATION

Leningradsky Prospect, 49-55, 125468 Moscow,

RUSSIAN FEDERATION

*E-mail address*: s.e.stepanov@mail.ru

DEPARTAMENT OF MATHEMATICS FINANCE UNIVERSITY

UNDER THE GAVERMENT OF RUSSIAN FEDERATION

Leningradsky Prospect, 49-55, 125468 Moscow,

RUSSIAN FEDERATION

*E-mail address*: i.i.tsyganok@mail.ru